\documentclass[12pt,leqno]{amsart}
\usepackage{amssymb,verbatim,enumerate,ifthen,cite}
\usepackage[mathscr]{eucal}
\usepackage[utf8]{inputenc}
\usepackage[T1]{fontenc}
\def\N{\mathbb{N}}
\def\R{\mathbb{R}}
\def\Q{\mathbb{Q}}
\def\Z{\mathbb{Z}}

\def\C{\mathbb{C}}

\newtheorem{theorem}{Theorem}[section]
\newtheorem*{theorem*}{Theorem}
\def\Thm#1#2{\ifthenelse{\equal{#1}{*}}{\begin{theorem*}#2\end{theorem*}}
	{\begin{theorem}\label{T#1}#2\end{theorem}}}

\def\thm#1{Theorem~\ref{T#1}}

\newtheorem{proposition}[theorem]{Proposition}
\newtheorem*{proposition*}{Proposition}
\def\Prp#1#2{\ifthenelse{\equal{#1}{*}}{\begin{proposition*}#2\end{proposition*}}
	{\begin{proposition}\label{P#1}#2\end{proposition}}}

\newtheorem{corollary}[theorem]{Corollary}
\newtheorem*{corollary*}{Corollary}
\def\Cor#1#2{\ifthenelse{\equal{#1}{*}}{\begin{corollary*}#2\end{corollary*}}
	{\begin{corollary}\label{C#1}#2\end{corollary}}}
\def\cor#1{Corollary~\ref{C#1}}

\newtheorem{lemma}[theorem]{Lemma}
\newtheorem*{lemma*}{Lemma}
\def\Lem#1#2{\ifthenelse{\equal{#1}{*}}{\begin{lemma*}#2\end{lemma*}}
	{\begin{lemma}\label{L#1}#2\end{lemma}}}

\theoremstyle{definition}
\newtheorem{remark}[theorem]{Remark}
\newtheorem*{remark*}{Remark}
\def\Rem#1#2{\ifthenelse{\equal{#1}{*}}{\begin{remark}\rm #2\end{remark}}
	{\begin{remark}\label{R#1}\rm #2\end{remark}}}

\newtheorem{statement}[theorem]{Statement}
\newtheorem*{statement*}{Statement}
\def\Stm#1#2{\ifthenelse{\equal{#1}{*}}{\begin{statement*}#2\end{statement*}}
	{\begin{statement}\label{s#1}#2\end{statement}}}

\newtheorem{example}[theorem]{Example}
\newtheorem*{example*}{Example}
\def\Exa#1#2{\ifthenelse{\equal{#1}{*}}{\begin{example*}#2\end{example*}}
	{\begin{example}\label{e#1} #2\end{example}}}

\def\eq#1{{\rm(\ref{E#1})}}
\def\Eq#1#2{\ifthenelse{\equal{#1}{*}}
	{\begin{equation*}\begin{aligned}[]#2\end{aligned}\end{equation*}}
	{\begin{equation}\begin{aligned}[]\label{E#1}#2\end{aligned}\end{equation}}}

\def\comment#1{}

\begin{document}
 \begin{flushright}
\end{flushright}
	\vspace{5mm}
	
	\date{\today}
	
\title{Cauchy--Schwarz-type inequalities for additive functions}
	
\author[Zs. P\'ales]{Zsolt P\'ales}
\address[Zs. P\'ales]{Institute of Mathematics, University of Debrecen, H-4002 Debrecen, Pf.\ 400, Hungary}
\email{pales@science.unideb.hu}
\author[M. K. Shihab]{Mahmood Kamil Shihab}
\address[M. K. Shihab]{Doctoral School of Mathematical and Computational Sciences, University of Debrecen, Hungary; Department of Mathematics, College of Education for Pure Sciences, University of Kirkuk, Iraq}
\email{mahmood.kamil@science.unideb.hu; mahmoodkamil30@uokirkuk.edu.iq}
	
\subjclass[2000]{Primary 39B22, 39B52, 39B62, 39B72}
\keywords{Levi--Civita-type functional equation, Cauchy--Schwarz-type functional inequality, additive function}
\thanks{The research of the first author was supported by the K-134191 NKFIH Grant.}

\begin{abstract}
The main goal of this paper is to show that if a real valued function defined on a groupoid satisfies a certain Levi--Civita-type functional equation, then it also fulfills a Cauchy--Schwarz-type functional inequality. 
In particular, if the groupoid is the multiplicative structure of commutative ring, then we can establish the existence of nontrivial additive functions satisfying inequalities connected to the multiplicative structure.
\end{abstract}
	
\maketitle
	
\section{Introduction}

In the theory of real and additive functions (see the monograph \cite{Kuc85} of Kuczma) there are several results which establish the existence of a discontinuous additive function which satisfies further algebraic conditions. One of the first problems of this kind was posed by Szabó \cite{Sza93} motivated by a question of Benz \cite{Ben90} and solved by Kominek, Reich, and Schwaiger \cite{KomReiSch98}. They proved that if $A:\R\to\R$ is an additive function which satisfies the equality $A(x)A(y)=0$ for all $(x,y)\in C$, where $C$ is a the unit circle, or is a hyperbola, or is an algebraic curve given by polynomials, then $A$ has to be equal to zero identically. Boros and Fechner \cite{BorFec15} and Boros, Fechner and Kutas \cite{BorFecKut16} extended these results to sets defined via generalized polynomials and to quadratic functions instead of additive ones, respectively, and they also examined the stability versions of such problems.

In \cite{BorFecKut16}, the case when $C$ is the graph of the hyperbola $xy=1$ was left open. Kanappan \cite[Chapter 1]{Kan09} proved that if, for some positive constant $a$, an additive function $A$ satisfies the condition $A(x)A(1/x)=a$ for all $x\neq0$, then $A$ has to be continuous. On the other hand, according to the remarks \cite{Ben80} and \cite{Ber81}, there exist discontinuous additive functions which fulfill the inequality $A(x)A(1/x)>0$ for all $x\neq 0$. On the other hand, using the theory of valuations of fields, Kutas \cite[Theorem 24]{Kut18} proved that there exists a nonzero (henceforth discontinuous) additive function which satisfies the equality $A(x)A(1/x)=0$ for all $x\neq 0$.

The above results motivated us to construct discontinuous additive real functions that enjoy properties that are connected to the multiplicative structure. It turned out that such properties could be possessed if the additive function satisfies Levi--Civita-type functional equations with respect to the multiplicative structure.

More generally, let $(G,*)$ be a groupoid. (Recall that a pair $(G,*)$ is said to be a groupoid if $\cdot$ is a binary operation on $G$, i.e., $\cdot:G\times G\to G$.) Let $A:G\to\R$ be a function such that there exist functions $f_1,\dots,f_n,g_1,\dots,g_n:G\to\R$ such that
the functional equation
\Eq{*}{
  A(x*y)=f_1(x)g_1(y)+\dots+f_n(x)g_n(y) \qquad(x,y\in G)
}
is fulfilled. Under certain assumptions on $n$ and on the functions $f_1,\dots,f_n,g_1,\dots,g_n$, we are going to prove that $A$ will satisfy either the inequality $A(x*y)^2\leq A(x*x)A(y*y)$ or the reversed one $A(x*x)A(y*y)\leq A(x*y)^2$. In the important particular case when the groupid is the multiplicative structure of a commutative ring and $A$ is additive, we will establish the existence of nontrivial additive functions which satisfy one of the above mentioned inequalities.

\section{The inequality $A(x*y)^2\leq A(x*x)A(y*y)$}

In our first result we assume that the function $A$ satisfies a Levi--Civita-type functional equation over a groupoid. 

\Thm{1+}{Let $(G,*)$ be a groupoid and let $A:G\to\R$ be a function. Assume that there exist $n\in\N$ and functions $f_1,\dots,f_n:G\to\R$ such that $A$ the Levi--Civita-type  functional equation
\Eq{1*=}{
  A(x*y)=f_1(x)f_1(y)+\dots+f_n(x)f_n(y)
}
holds for all $x,y\in G$. Then, $A$ satisfies the functional inequality 
\Eq{1*<}{
  A(x*y)^2\leq A(x*x)A(y*y)
}
for all $x,y\in G$.}

\begin{proof} Let $x,y\in G$. In view of the functional equation \eq{1*=}, the inequality \eq{1*<} can be rewritten as
\Eq{*}{
  \big(f_1(x)f_1(y)+\dots+f_n(x)f_n(y)\big)^2
  \leq \big(f_1(x)^2+\dots+f_n(x)^2\big)\big(f_1(y)^2+\dots+f_n(y)^2\big),
}
which follows from the Cauchy--Schwarz inequality when we apply it to the $n$-dimensional vectors $(f_1(x),\dots,f_n(x))$ and $(f_1(y),\dots,f_n(y))$.
\end{proof}

If the groupoid is the multiplicative semigroup of a commutative ring $(R,+,\cdot)$ and $A$ is additive, then we can establish a characterization of the corresponding inequality. Recall that in a ring, the product $x\cdot y$ of the elements $x,y\in R$ is simply denoted by $xy$, and $x^2$ is defined to be the product $x\cdot x$. 

\Thm{1}{Let $(R,+,\cdot)$ be a commutative ring and let $A:R\to\R$ be an additive function. If one of the following conditions 
\begin{enumerate}[(i)]
 \item $A=0$,
 \item $A(x^2)\geq0$ for all $x\in R$,
 \item $A(x^2)\leq0$ for all $x\in R$
\end{enumerate}
hold, then $A$ satisfies the inequality 
\Eq{1<}{
  A(xy)^2\leq A(x^2)A(y^2)
}
for all $x,y\in R$. Conversely, if $R$ has a multiplicative unit element and $A$ satisfies the inequality \eq{1<} for all $x,y\in R$, then one of the conditions (i), (ii), and (iii) must be satisfied.}

\begin{proof} 
The inequality \eq{1<} is obvious if $A=0$, i.e., if condition (i) holds.

Now assume that $A$ satisfies condition (ii) and let $x,y\in R$ be fixed. Then, for all $n\in\N$ and $k\in\Z$, we get that
\Eq{*}{
  0\leq A((nx+ky)^2)
  =A(n^2x^2+2nkxy+k^2y^2)
  =n^2A(x^2)+2nkA(xy)+k^2A(y^2).
}
Dividing this inequality by $n^2$, we can conclude that
\Eq{*}{
  0\leq A(x^2)+2\frac{k}{n}A(xy)+\frac{k^2}{n^2}A(y^2).
}
Because $n\in\N$ and $k\in\Z$ were arbitrary, we obtain that 
\Eq{*}{
  0\leq A(x^2)+2rA(xy)+r^2A(y^2)
}
is valid for all rational number $r$. By the density of rational numbers, it follows that the above inequality is true for all real number $r$. The polynomial on the right hand side cannot have two distinct real roots, therefore, its discriminant has to be non positive, i.e.,
\Eq{*}{
  (2A(xy))^2-4A(x^2)A(y^2)\leq0.
}
This inequality reduces to \eq{1<}.

In the case when condition (iii) holds, then the additive function $(-A)$ satisfies condition (ii) and hence the inequality \eq{1<} holds with $(-A)$ instead of $A$, which again shows that \eq{1<} is valid.

To verify the reversed implication, assume that $R$ possesses a multiplicative unit element which will be denoted by $e$ and assume that $A$ satisfies inequality \eq{1<} for all $x,y\in R$. Substituting $y:=e$ in \eq{1<}, it follows that 
\Eq{1x}{
0\leq A(x)^2\leq A(x^2)A(e)
} 
holds for all $x\in R$. We can now distinguish three cases according to the possibilities $A(e)=0$, $A(e)>0$, and $A(e)<0$.

If $A(e)=0$, then \eq{1x} yields that $A(x)=0$ for all $x\in R$, i.e., condition (i) is valid.

In the case when $A(e)>0$, it follows that $A$ satisfies condition (ii), while in the case when $A(e)<0$, we can see that $A$ satisfies condition (iii).
\end{proof}

\section{The inequality $A(x*x)A(y*y)\leq A(x*y)^2$} 

In the subsequent two theorems, we present to Levi--Civita-type functional equations which imply the inequality in the title of this section.

\Thm{3}{Let $(G,*)$ be a groupoid and $A:G\to\R$ be a function. Assume that there exist two functions $f,g:G\to\R$ such that the Levi--Civita-type functional equation
\Eq{3=}{
  A(x*y)=f(x)f(y)-g(x)g(y)
}
holds for all $x,y\in G$. Then $A$ satisfies the functional inequality 
\Eq{3<}{
 A(x*x)A(y*y)\leq A(x*y)^2
}
for all $x,y\in G$.}

\begin{proof} Let $x,y\in G$. According to the functional equation \eq{3=}, the inequality \eq{3<} can be rewritten as
\Eq{*}{
  (f(x)^2-g(x)^2)(f(y)^2-g(y)^2)
  \leq (f(x)f(y)-g(x)g(y))^2.
}
Observe that this inequality is equivalent to 
\Eq{*}{
 0\leq(g(x)f(y)-f(x)g(y))^2,
}
which is obviously valid.
\end{proof}

\Thm{4}{Let $(G,*)$ be a groupoid. Let $A:R\to\R$ be a function. Assume that there exist $f,g:R\to\R$ such that the Levi--Civita-type functional equation
\Eq{4=}{
  A(x*y)=f(x)g(y)+g(x)f(y)
}
holds for all $x,y\in G$. Then, for all $x,y\in G$, $A$ satisfies the functional inequality \eq{3<}.}

\begin{proof} Let $x,y\in G$. According to the functional equation \eq{4=}, the inequality \eq{3<} can be rewritten as
\Eq{*}{
  4f(x)g(x)f(y)g(y)\leq (f(x)g(y)+g(x)f(y))^2.
}
Observe that this inequality is equivalent to 
\Eq{*}{
 0\leq(g(x)f(y)-f(x)g(y))^2,
}
which is obviously valid.
\end{proof}

\Cor{4}{Assume that $A:\R\to\R$ satisfies the Leibniz Rule, i.e.,
\Eq{*}{
  A(xy)=xA(y)+A(x)y \qquad (x,y\in\R).
}
Then, for all $x,y\in\R$, the inequality
\Eq{C4<}{
A(x^2)A(y^2)\leq A(xy)^2
}
holds.}

\begin{proof}
Observe that with groupoid $(G,*):=(\R,\cdot)$ and with the notations $g:=A$ and $f(x):=x$, $(x\in G)$, the equality \eq{4=} of \thm{4} holds. Therefore, $A$ satisfies inequality \eq{3<} for all $x,y\in\R$, hence \eq{C4<} is also satisfied.
\end{proof}

In particular, if $A:\R\to\R$ is a derivation (i.e., $A$ is additive and satisfies the Leibniz Rule with respect, then the above corollary implies that it fulfills the inequality \eq{C4<}. 

If the groupoid is the multiplicative semigroup of a commutative ring $(R,+,\cdot)$ and $A$ is additive, then we can establish a characterization of the inequality \eq{3<} over a particular subset of the ring. 

\Thm{2}{Let $(R,+,\cdot)$ be a commutative ring with a multiplicative unit element $e$ and $A:R\to\R$ be an additive function with $A(e)\neq0$. Let the subset $R_A\subseteq R$ be defined by
\Eq{*}{
  R_A:=\{x\in R\mid  0\leq A(x^2)A(e)\}
}
Then $e\in R_A$ and $A$ satisfies the following functional inequality
\Eq{2<}{
  A(x^2)A(y^2)\leq A(xy)^2
}
for all $x,y\in R_A$ if and only if
\Eq{2xe}{
  A(x^2)A(e)\leq A(x)^2
}
for all $x\in R_A$.}

\begin{proof} The inclusion $e\in R_A$ is obvious. Now, putting $y:=e$, we can see that the inequality \eq{2<} implies \eq{2xe}.

To prove the reversed implication, assume that \eq{2xe} is valid for all $x\in R_A$. Then it is also valid for all $x\in R$, since, for $x\in R\setminus R_A$, the left hand side of the inequality is negative, while the right hand side is nonnegative. Introduce the function $A_0:=A/A(e)$. Then, $A_0$ is additive, $A_0(e)=1$ and, dividing \eq{2xe} by $A(e)^2>0$ side by side, for all $x\in R$, we get that
\Eq{2x1}{
  A_0(x^2)\leq A_0(x)^2
}
Let $x,y\in R_A$ be fixed and $n\in\N$, $k\in\Z$ be arbitrary. Then, \eq{2x1} yields 
\Eq{*}{
    A_0((nx+ky)^2)\leq A_0(nx+ky)^2.
}
Using the additivity of $A_0$, we get
\Eq{*}{
    n^2A_0(x^2)+2nkA_0(xy)+k^2A_0(y^2)
    \leq n^2A_0(x)^2+2nkA_0(x)A_0(y)+k^2A_0(y)^2.
}
Dividing this inequality by $n^2$, we obtain
\Eq{*}{
    A_0(x^2)+2\tfrac{k}{n}A_0(xy)+\big(\tfrac{k}{n}\big)^2A_0(y^2)
    \leq A_0(x)^2+2\tfrac{k}{n}A_0(x)A_0(y)+\big(\tfrac{k}{n}\big)^2A_0(y)^2.
}
Therefore, for any rational number $r\in\Q$, 
\Eq{*}{
    0\leq (A_0(x)^2-A_0(x^2))+2r(A_0(x)A_0(y)-A_0(xy))+r^2(A_0(y)^2-A_0(y^2)).
}
Using the continuity of both sides as a function of $r$, it follows that the same inequality is valid for all $r\in\R$. 
Thus, the discriminant of this quadratic polynomial has to be nonpositive, i.e.,
\Eq{A9}{
  (A_0(x)A_0(y)-A_0(xy))^2\leq (A_0(x)^2-A_0(x^2))(A_0(y)^2-A_0(y^2))
}
and hence
\Eq{*}{
  |A_0(x)A_0(y)-A_0(xy)|\leq \sqrt{(A_0(x)^2-A_0(x^2))(A_0(y)^2-A_0(y^2))}=Q(x)Q(y),
}
thus
\Eq{*}{
  \big||A_0(x)A_0(y)|-|A_0(xy)|\big|\leq
  |A_0(x)A_0(y)-A_0(xy)|\leq Q(x)Q(y),
}
where $Q(u):=\sqrt{A_0(u)^2-A_0(u^2)}\geq0$ ($u\in R$). Then, for all $u\in R$,
\Eq{A7}{
	A_0(u)^2=Q(u)^2+A_0(u^2).
}
Therefore $|A_0(xy)|$ satisfies the inequality
\Eq{A0}{
  |A_0(x)A_0(y)|-Q(x)Q(y)\leq |A_0(xy)| \leq
  |A_0(x)A_0(y)|+Q(x)Q(y).
}
We are going to show that
\Eq{A8}{
  A_0(x^2)A_0(y^2)\leq A_0(xy)^2.
}

To see this inequality, we will prove that 
\Eq{A2}{
	Q(x)^2Q(y)^2\leq A_0(x)^2A_0(y)^2.
}
and 
\Eq{A3}{
  A_0(x^2)A_0(y^2)\leq\Big(|A_0(x)A_0(y)|-Q(x)Q(y)\Big)^2.
}
Since $x$ and $y$ belong to $R_A$, therefore, we have that $A_0(x^2)\geq0$ and $A_0(y^2)\geq0$, then, with $u\in\{x,y\}$, the equality \eq{A7} implies that
\Eq{*}{
  Q(x)^2\leq A_0(x)^2 \qquad\mbox{and}\qquad
  Q(y)^2\leq A_0(y)^2.
}
Multiplying these inequalities side by side, we get that \eq{A2} holds. Therefore, we can conclude that
\Eq{A4}{
  Q(x)Q(y)\leq |A_0(x)A_0(y)|
}
which is equivalent to \eq{A2}.

By the obvious inequality
\Eq{*}{
  \big(|A_0(x)|Q(y)-|A_0(y)|Q(x)\big)^2\geq 0,
}
we have that 
\Eq{A6}{
  2|A_0(x)A_0(y)|Q(x)Q(y)\leq A_0(x)^2Q(y)^2+A_0(y)^2Q(x)^2.
}
Therefore using the equality \eq{A7} with $u\in\{x,y\}$ and the inequality \eq{A6}, we obtain
\Eq{*}{
  A_0(x^2)A_0(y^2)&= A_0(x)^2A_0(y)^2-A_0(x)^2Q(y)^2-Q(x)^2A_0(y)^2+Q(x)^2Q(y)^2\\
  &\leq A_0(x)^2A_0(y)^2-2|A_0(x)A_0(y)|Q(x)Q(y)+Q(x)^2Q(y)^2\\
  &=(|A_0(x)A_0(y)|-Q(x)Q(y))^2.
}
This shows that the inequality \eq{A3} holds.

In view of \eq{A4}, the first inequality in \eq{A0} implies that
\Eq{*}{
	(|A_0(x)A_0(y)|-Q(x)Q(y))^2\leq A_0(xy)^2.
}
This, combined with the inequality \eq{A3} yields that \eq{A8} is valid, indeed. Therefore,
\Eq{*}{
  A(x^2)A(y^2)=A_0(x^2)A_0(y^2)A(e)^2\leq A_0^2(xy)A(e)^2=A(xy)^2,
}
which completes the proof of the inequality \eq{2<} for $x,y\in R_A$.
\end{proof}

In the following example we show that the additivity of the function $A$ in \thm{2} is necessary.
\Exa{a}{
Let $q\in(0,1)$ and let $f:\R\to\R$ be a function defined by 
 \Eq{*}{
 	f(x)=\begin{cases}
 		 x & x\neq 1,\\
 		 q & x=1.
 		 \end{cases}
 }
Clearly, $f$ is not additive. Therefore, for $x\not\in\{1,-1\}$, we have that 
\Eq{*}{
 f(x^2)f(1)=qx^{2}\leq x^2=f(x)^2
}
for $x=\pm 1$, 
\Eq{*}{
 f(1^2)f(1)=q^2=f(1)^2, \qquad
 f((-1)^2)f(1)=q^2<1=f(-1)^2,
}
which shows that \eq{2xe} is satisfied for all $x\in\R$.
On the other hand, for $x,y\in\R\setminus\{1,-1\}$ with $xy=1$, we can conclude that 
\Eq{*}{
f(x^2)f(y^2)=x^2y^2=1>q^2=f(xy)^2,
}
which shows that \eq{2<} is not satisfied.}

The next example shows that if the function $A$ in \thm{2} is non-additive, continuous and satisfies $A(e)=0$, then the conclusion of \thm{2} may not be valid. 

\Exa{co}{
Let $A:\R\to\R$ be defined as $A(x)=|x-1|$. Note that $A$ is continuous and not additive. Since $A(1)=0$ this implies that
\Eq{*}{
  A(x^2)A(1)=0\leq (x-1)^2=A(x)^2.
}
Thus the inequality \eq{2xe} holds for all $x\in\R$. On the other hand we have that
\Eq{*}{
  A(x^2)A(y^2)=|x^2-1||y^2-1|\qquad\mbox{and}\qquad A(xy)^2=(xy-1)^2.
}
Hence for $x=2$ and $y=\frac{1}{2}$ we have that $A(x^2)A(y^2)=\frac{9}{4}$ but $A(xy)^2=0$, therefore the inequality \eq{2<} does not hold.
}

\section{Consequences of systems of Levi--Civita-type functional equations}

\Thm{5}{Let $(G,*)$ be a groupoid and $A,B:G\to\R$ be  functions. Assume that there exist $f,g:G\to\R$ such that $A$ and $B$ satisfy the following system of Levi--Civita-type functional equations
\Eq{2=}{
  A(x*y)&=f(x)f(y)-g(x)g(y) \qquad\mbox{and}\\
  B(x*y)&=f(x)g(y)+g(x)f(y)  
}
for all $x,y\in G$. Then the inequalities
\Eq{xy1}{
  -B(x*y)^2\leq A(x*x)A(y*y)\leq A(x*y)^2
}
and 
\Eq{xy2}{
  -A(x*y)^2\leq B(x*x)B(y*y)\leq B(x*y)^2
}
hold for all $x,y\in G$.}

\begin{proof}
In view of two functional equations in \eq{2=}, for $x,y\in G$, we have that
\Eq{*}{
B(x*y)^2&+A(x*x)A(y*y)\\
&=f(x)^2g(y)^2+2f(x)g(y)g(x)f(y)+g(x)^2f(y)^2
+(f(x)^2-g(x)^2)(f(y)^2-g(y)^2)\\
&=(f(x)f(y)+g(x)g(y))^2\geq0,
}
which proves the left hand side inequality in \eq{xy1}.
The right hand side inequality in \eq{xy1} is a direct consequence of \thm{3}.

Again, in view of two equations in \eq{2=}, for $x,y\in G$ we have that
\Eq{*}{
A(x*y)^2&+B(x*x)B(y*y)\\
&=f(x)^2f(y)^2-2f(x)f(y)g(x)g(y)+g(x)^2g(y)^2+4f(x)g(x)f(y)g(y)\\
&=(f(x)f(y)+g(x)g(y))^2\geq 0.
}
This implies the left hand side inequality in \eq{xy2}.
On the other hand, applying \thm{4} for the function $B$ instead of $A$, we obtain that
\Eq{*}{
	B(x*x)B(y*y)\leq B(x*y)^2. 
}
This shows that the second inequality of \eq{xy2} holds for $x,y\in G$.
\end{proof}

An interesting consequence of the functional equations in \eq{2=} is that $A$ and $B$ satisfy the following identity:
\Eq{*}{
  B(x*y)^2+A(x*x)A(y*y)=A(x*y)^2+B(x*x)B(y*y) \qquad(x,y\in G).
}
Therefore, the inequalities \eq{xy1} and \eq{xy2} can be expressed as the following chain of inequalities
\Eq{*}{
  0&\leq A(x*x)A(y*y)+(B(x*y))^2\\
  &=B(x*x)B(y*y)+A(x*y)^2
  \leq A(x*y)^2+B(x*y)^2 \qquad(x,y\in G).
}

\Cor{5}{For all $x,y\in\R$, we have 
\Eq{*}{
  -\sin(x+y)^2&\leq\cos(2x)\cos(2y)\leq\cos(x+y)^2
  \qquad\mbox{and}\\
  -\cos(x+y)^2&\leq\sin(2x)\sin(2y)\leq\sin(x+y)^2.
}}

\begin{proof}
Observe that the trigonometric functions $\cos:\R\to\R$ and $\sin:\R\to\R$ satisfy the functional equations  
\Eq{*}{
  \cos(x+y)&=\cos(x)\cos(y)-\sin(x)\sin(y) \qquad\mbox{and}\\
  \sin(x+y)&=\sin(x)\cos(y)+\cos(x)\sin(y)  
}
for all $x,y\in\R$. Therefore, \eq{2=} holds with $A:=f:=\cos$ and $B:=g:=\sin$ over the groupoid $(\R,+)$. Consequently, \eq{xy1} and \eq{xy2} are satisfied for all $x,y\in\R$, which imply the assertion.
\end{proof}

\Cor{A1}{Let $(G,*)$ be a groupoid and let $\varphi:G\to\C$ be a homomorphism into the multiplicative semigroup of complex numbers. Define $A:=\Re\varphi$ and $B:=\Im\varphi$. Then, for all $x,y\in G$, the inequalities \eq{xy1} and \eq{xy2} hold.}

\begin{proof} 
Using the multiplicativity of $\varphi$,
for all $x,y\in G$, we get that
\Eq{*}{
  A(x*y)&=\Re(\varphi(x*y))
  =\Re(\varphi(x)\varphi(y))\\
  &=\Re((A(x)+iB(x))(A(y)+iB(y)))
  =A(x)A(y)-B(x)B(y),\\
  B(x*y)&=\Im(\varphi(x*y))
  =\Im(\varphi(x)\varphi(y))\\
  &=\Im((A(x)+iB(x))(A(y)+iB(y)))
  =A(x)B(y)+B(x)A(y).
}
Therefore, the functional equations in \eq{2=} are satisfied with $f:=A$ and $g:=B$. Thus, according to \thm{5}, we obtain that the inequalities \eq{xy1} and \eq{xy2} hold for all $x,y\in G$, which was to be shown.
\end{proof}

\Cor{A2}{Let $\varphi:\C\to\C$ be an automorphism of the field $\C$.  Define $A:=\Re\varphi$ and $B:=\Im\varphi$. Then $A:\C\to\R$ and $B:\C\to\R$ are additive mappings, furthermore, for all $x,y\in\C$,
\Eq{xy+}{
  -B(xy)^2\leq A(x^2)A(y^2)\leq A(xy)^2
  \quad\mbox{and}\quad
  -A(xy)^2\leq B(x^2)B(y^2)\leq B(xy)^2.
}}

\begin{proof} Using the additivity of $\varphi$, for all $x,y\in\C$, we obtain that
\Eq{*}{
  A(x+y)&=\Re(\varphi(x+y))
  =\Re(\varphi(x)+\varphi(y))\\
  &=\Re((A(x)+iB(x))+(A(y)+iB(y)))
  =A(x)+A(y),\\
  B(x+y)&=\Im(\varphi(x+y))
  =\Im(\varphi(x)+\varphi(y))\\
  &=\Im((A(x)+iB(x))+(A(y)+iB(y)))
  =B(x)+B(y).
}
These equalities show that $A:\C\to\R$ and $B:\C\to\R$ are additive mappings. 

By the multiplicativity of $\varphi$, it maps the groupoid $(G,*):=(\C,\cdot)$ into itself. Thus, according to  \cor{A1}, we obtain that the inequalities \eq{xy1} and \eq{xy2} hold for all $x,y\in\C$. This yields the assertion.
\end{proof}

\comment{Furthermore, in the case when $\varphi(i)=i$, for all $x\in\C$,
\Eq{*}{
	B(ix)=\Im\varphi(ix)=\Im(\varphi(i)\varphi(x))
	=\Im(i\varphi(x))=\Im(iA(x)-B(x))=A(x).
}}

The following result is a counterpart of \thm{5}.

\Thm{6}{Let $(G,*)$ be a groupoid and $A,B:G\to\R$ be  functions. Assume that there exist $f,g:G\to\R$ such that $A$ and $B$ satisfy the following Levi--Civita-type functional equations
\Eq{6=}{
  A(x*y)&=f(x)f(y)+g(x)g(y) \qquad\mbox{and}\\
  B(x*y)&=f(x)g(y)+g(x)f(y)  
}
for all $x,y\in G$. Then the inequalities
\Eq{6xy1}{
  B(x*x)B(y*y)\leq A(x*y)^2\leq A(x*x)A(y*y)
}
and 
\Eq{6xy2}{
  B(x*x)B(y*y)\leq B(x*y)^2\leq A(x*x)A(y*y).
}
hold for all $x,y\in G$.}

\begin{proof}
In view of two functional equations in \eq{6=}, for $x,y\in G$, we have that
\Eq{*}{
A(x*y)^2&-B(x*x)B(y*y) \\
&=f(x)^2f(y)^2+2f(x)f(y)g(x)g(y)+g(x)^2g(y)^2-4f(x)g(x)f(y)g(y)\\
&=(f(x)f(y)-g(x)g(y))^2\geq 0.
}
This implies the left hand side inequality in \eq{6xy1}.
On the other hand, applying \thm{1+} for the function $A$ and $n=2$, $f_1:=f$, $f_2:=g$, we obtain that the second inequality of \eq{6xy1} holds for $x,y\in G$.

Again, in view of two equations in \eq{6=}, for $x,y\in G$, we have that
\Eq{*}{
A(x*x)&A(y*y)-B(x*y)^2\\
&=(f(x)^2+g(x)^2)(f(y)^2+g(y)^2)-f(x)^2g(y)^2-2f(x)g(y)g(x)f(y)-g(x)^2f(y)^2\\
&=(f(x)f(y)-g(x)g(y))^2\geq0,
}
which proves the right hand side inequality in \eq{6xy2}.
The left hand side inequality in \eq{6xy2} is a direct consequence of \thm{4} (applied to $B$ instead of $A$).
\end{proof}

An interesting consequence of the functional equations in \eq{2=} is that $A$ and $B$ satisfy the following identity:
\Eq{*}{
  B(x*x)B(y*y)+A(x*x)A(y*y)=A(x*y)^2+B(x*y)^2 \qquad(x,y\in G).
}

\Cor{6}{For all $x,y\in\R$, we have 
\Eq{C6}{
  \sinh(2x)\sinh(2y)
  \leq\sinh(x+y)^2
  <\cosh(x+y)^2
  \leq\cosh(2x)\cosh(2y)
}}

\begin{proof}

Observe that the hyperbolic functions $\cosh:\R\to\R$ and $\sinh:\R\to\R$ satisfy the functional equations  
\Eq{*}{
  \cosh(x+y)&=\cosh(x)\cosh(y)+\sinh(x)\sinh(y) \qquad\mbox{and}\\
  \sinh(x+y)&=\sinh(x)\cosh(y)+\cosh(x)\sinh(y)  
}
for all $x,y\in\R$. Therefore, \eq{6=} holds with $A:=f:=\cosh$ and $B:=g:=\sinh$ over the groupoid $(\R,+)$. Consequently, \eq{6xy1} and \eq{6xy2} are satisfied for all $x,y\in\R$, which imply the first and last inequalities in \eq{C6}.
The central inequality follows from the identity
$\cosh^2-\sinh^2=1$.
\end{proof}

To formulate the next result, let $p$ be a square free positive integer and let $\Q(\sqrt{p})$ denote the subfield of $\R$ generated by $\sqrt{p}$. Then, one can see that $\Q(\sqrt{p})=\{a+b\sqrt{p}: a,b\in\Q\}$.

\Thm{B}{Let $p$ be a square free positive integer. Then there exist two discontinuous additive functions  $A:\Q(\sqrt{p})\to\R$ and $B:\Q(\sqrt{p})\to\R$ such that the functions $A$ and $B$ fulfill the inequalities
	\Eq{a}{
		B(x^2)B(y^2)&\leq A(xy)^2
		\leq A(x^2)A(y^2),\\
		B(x^2)B(y^2)&\leq B(xy)^2
		\leq A(x^2)A(y^2) 
	}
for all $x,y\in\Q(\sqrt{p})$.
}
\begin{proof}
Define the functions $A:\Q(\sqrt{p})\to\R$ and $B:\Q(\sqrt{p})\to\R$ by
	\Eq{*}{
		A(a+b\sqrt{p}):=a\qquad \mbox{and}\qquad B(a+b\sqrt{p}):=b\sqrt{p} \qquad(a,b\in\Q).
	}
	
We show that $A$ and $B$ are additive. Indeed, let $x=a_1+b_1\sqrt{p}$ and $y=a_2+b_2\sqrt{p}$ be two arbitrary points of $\Q(\sqrt{p})$, where $a_1,a_2,b_1,b_2\in\Q$. According to the definition of $A$ and $B$, we have that
	\Eq{*}{
		A(x+y)&=A(a_1+a_2+(b_1+b_2)\sqrt{p})
		=a_1+a_2=A(x)+A(y)
	}
	and
	\Eq{*}{
		B(x+y)&=B(a_1+a_2+(b_1+b_2)\sqrt{p})
		=(b_1+b_2)\sqrt{p}=B(x)+B(y).
	}
	This proves that $A$ and $B$ are additive, indeed. 

Next we prove that $A$ and $B$ satisfy the functional equations in \eq{6=} with $f:=A$, and $g:=B$, where the groupoid $(G,*)$ is equal to $(\Q(\sqrt{p}),\cdot)$. Indeed, let $x=a_1+b_1\sqrt{p}$ and $y=a_2+b_2\sqrt{p}$ be two arbitrary points of $\Q(\sqrt{p})$, where $a_1,a_2,b_1,b_2\in\Q$. Then
	\Eq{*}{
	 A(xy)=A((a_1a_2+b_1b_2p)+(a_1b_2+a_2b_1)\sqrt{p})
	 =a_1a_2+b_1b_2p=A(x)A(y)+B(x)B(y)
	}
and similarly,
	\Eq{*}{
	 B(xy)=B((a_1a_2+pb_1b_2)+(a_1b_2+a_2b_1)\sqrt{p})
	 =(a_1b_2+a_2b_1)\sqrt{p}=A(x)B(y)+B(x)A(y).
	}
Therefore, according to \thm{6}, the inequalities \eq{6xy1} and \eq{6xy2} holds, which prove that the inequalities in \eq{a} are also valid.

Finally, we show that the functions $A$ and $B$ are discontinuous at the point $u:=1+\sqrt{p}$. By the density of the set $\Q$ in $\R$, there exists a sequence $(x_n)$ of rational numbers converging to $u/\sqrt{p}$. Then the sequence $(x_n\sqrt{p})$ converges to $u$. We have that $A(u)=1$ but, for all $n\in\N$, $A(x_n\sqrt{p})=0$, therefore $A$ is discontinuous at $u$. Furthermore, there exists a sequence $(y_n)$ of rational numbers, which converges to $u$. Since $B(u)=\sqrt{p}$ and $B(y_n)=0$ for all $n\in\N$, thus we can conclude that the function $B$ is also discontinuous at $u$.
\end{proof}


\providecommand{\bysame}{\leavevmode\hbox to3em{\hrulefill}\thinspace}
\providecommand{\MR}{\relax\ifhmode\unskip\space\fi MR }
\providecommand{\MRhref}[2]{%
  \href{http://www.ams.org/mathscinet-getitem?mr=#1}{#2}
}
\providecommand{\href}[2]{#2}

\end{document}